\documentclass[11pt]{amsart}
\usepackage{amsmath,amsfonts,amssymb}
\newcommand{\e}{\mathrm{e}}
\newcommand{\esssup}{\mathop{\mathrm{ess\,sup}}}
\newcommand{\rr}{\mathbb{R}}
\newtheorem{theorem}{Theorem}{\rm}
\newtheorem{lemma}[theorem]{Lemma}{\rm}
\newtheorem{prop}[theorem]{Proposition}{\rm}
\newtheorem{cor}[theorem]{Corollary}{\rm}
\newtheorem{ex}{Example}
\newtheorem{rem}{Remark}
\numberwithin{equation}{section}

\def\id{{\rm ri}\,\mathcal{D}}

\def\R{\mathbb{R}}

\def\N{\mathbb{N}}

\def\Q{\mathbf{Q}}

\def\P{\mathbf{P}}
\def\Z{\mathbf{Z}}

\def\A{\mathbf{A}}

\def\B{\mathbf{B}}
\def\c{\mathbf{C}}

\def\X{\mathbf{X}}

\def\z{\mathbf{z}}

\def\bd{\mathbf{d}}
\def\blambda{\mathbf{\lambda}}
\def\x{\mathbf{x}}
\def\y{\mathbf{y}}

\def\b{\mathbf{b}}
\def\c{\mathbf{c}}

\def\rix{{\rm int}\,\X}

\begin{document}

\title[$L^p$-norms and Cramer transform  in optimization]
{$L^p$-norms, log-barriers and Cramer transform  in optimization}

\author{J.B. Lasserre}

\address{LAAS-CNRS and Institute of Mathematics\\
University of Toulouse\\
LAAS, 7 avenue du Colonel Roche\\
31077 Toulouse C\'edex 4,France}
\email{lasserre@laas.fr}
\author{E.S. Zeron}
\address{Depto. Matematicas\\
CINVESTAV-IPN, Apdo. Postal 14-740\\
Mexico, D.F. 07000, Mexico}
\email{eszeron@math.cinvestav.mx}

\begin{abstract}
We show that the Laplace approximation of a supremum by $L^p$-norms
has interesting consequences in optimization. For instance, the logarithmic barrier functions (LBF)
of a primal convex problem $\P$ and its dual $\P^*$ appear naturally
when using this simple approximation technique for the value function $g$ of $\P$
or its Legendre-Fenchel conjugate $g^*$. In addition,
minimizing the LBF of the dual $\P^*$ is just evaluating 
the Cramer transform of the Laplace approximation of $g$. Finally,
this technique permits to sometimes define an explicit dual 
problem $\P^*$ in cases when the Legendre-Fenchel conjugate $g^*$ cannot be 
derived explicitly from its definition.
\end{abstract}

\keywords{optimization; Logarithmic Barrier Function; Legendre-Fenchel
and Cramer transforms}

\subjclass{90C25 90C26}\maketitle

\section{Introduction}

Let $f:\X\to\rr$ and $\omega:\X\to\rr^m$ be a pair of continuous 
mappings defined on the convex cone $\X\subseteq\rr^n$.
Consider the function $g:\rr^m\to\rr\cup\{-\infty\}$ given by the formula:
\begin{equation}\label{fenchel1}
\y\mapsto g(\y)\,:=\,\sup_\x\:\{f(\x):\omega(\x)\leq\y,\,\x\in \X\}.
\end{equation}
For each fixed $\y\in\rr^m$, computing $g(\y)$ 
is solving the optimization problem
\begin{equation}
\label{pbprimal}
\P:\quad \sup_\x\:\{f(\x):\omega(\x)\leq\y,\,\x\in \X\},\end{equation}
and $g$ is called the {\it value} function associated with $\P$. The 
value function $g$ provides a systematic way to generate a dual problem 
$\P^*$ via its Legendre-Fenchel conjugate denoted $g^*:\R^m\to\R\cup\{-\infty\}$.
In the
concave version (i.e. when $g$ and $g^*$ are concave instead of convex), the 
Legendre-Fenchel conjugate $g^*$ is defined by
\begin{equation}\label{lfg}
\lambda\mapsto g^*(\lambda)\,:=\,\inf_{\y\in\rr^m}\,\{\lambda'\y-g(\y)\:\},
\end{equation}
and is finite on some domain $\mathcal{D}\subset\R^m_+$.
Then a dual problem is defined by:
\begin{equation}
\label{dual1}
\P^*:\:\tilde{g}(\y)\,:=\,(g^*)^*(\y)\,=\,\inf_\lambda\:\{ \lambda'\y-g^*(\lambda)\}.
\end{equation}
Of course, one has the property $\tilde{g}(\y)\geq g(\y)$
because from
\[g^*(\lambda)\,=\,\inf_{\x} \{\lambda'\x-g(\x)\,\}\,\leq\,\lambda'\y-g(\y),\]
one may deduce
\[\tilde{g}(\y)\,=\,(g^*)^*(\y)\,=\,\inf_\lambda \{\lambda'\y-g^*(\lambda)\,\}\,\geq\,
\inf_\lambda\,\{\lambda'\y+g(\y)-\lambda'\y\}\,=\,g(\y).\]
Moreover, notice that:
\begin{eqnarray}
\nonumber
\tilde{g}(\y)
&=&\inf_\lambda\:\{ \lambda'\y-g^*(\lambda)\}\,=\,\inf_\lambda\:\{ \lambda'\y+\sup_\z \{g(\z)-\lambda'\z\}\:\}\\
\nonumber&=&\inf_\lambda\:\{ \lambda'\y+\sup_\z\,\{ \sup_{\x\in\X} \{f(\x)\::\omega(\x)\leq \z\}-\lambda'\z\}\:\}\\
\nonumber&=&\left\{\begin{array}{l}
+\infty\quad\mbox{if $\lambda\not\in\R^m_+$}\\
\displaystyle\inf_{\lambda\in\R^m_+}\:\{ \lambda'\y+ \sup_{\x\in\X} \{f(\x)-\lambda'\omega(\x)\}\:\}\:\mbox{otherwise}\end{array}\right.\\
\label{dual2}&=&\inf_{\lambda\in\R^m_+}\:\sup_{\x\in\X}\:\{ f(\x)+\lambda (\y-\omega(\x)\,\}.
\end{eqnarray}
Finally, noting that $\displaystyle\inf_{\lambda\in\R^m_+}\:\{ f(\x)+\lambda' (\y-\omega(\x)\,\}=f(\x)$ if $\omega(\x)\leq \y$ and $-\infty$ otherwise, one may write
\begin{eqnarray}
\label{min-max1}
g(\y)&=&\sup_{\x\in\X}\:\inf_{\lambda\in\R^m_+}\:\{ f(\x)+\lambda' (\y-\omega(\x)\,\}\\
\label{min-max2}&\leq&\tilde{g}(\y)=\inf_{\lambda\in\R^m_+}\:\sup_{\x\in\X}\:\{ f(\x)+\lambda' (\y-\omega(\x)\,\}\quad\mbox{[by (\ref{dual2})]},
\end{eqnarray}
and the equality $g(\y)=\tilde{g}(\y)$ holds true
under some convexity assumption. 

However, in general $g^*$ cannot be obtained {\it explicitly} from its definition (\ref{lfg}), 
and for dual methods to solve $\P$,
the inner maximization in (\ref{min-max2})
must be done numerically for each fixed $\lambda$. A notable exception is the conic optimization 
problem where $f$ and $\omega$ are both linear mappings, for which 
the dual (\ref{dual1}) has an explicit form in terms of $\lambda$. Of course, alternative explicit duals 
have been proposed but they involve both primal ($\x$) and dual 
($\blambda$) variables. In particular, the Wolfe \cite{wolfe} and 
Mond-Weir \cite{mond1} duals even allow to consider weakened 
notions of convexity like e.g. pseudo- or quasi-convexity. For a nice 
exposition and related references on this topic, the interested reader 
is referred to Mond \cite{mond2} and the references therein.

\smallskip
{\bf Contribution.} Our contribution is to show
that the simple and well-known Laplace approximation of a supremum via a
converging sequence of $L^p$-norms has interesting consequences in optimization,
for both primal and dual problems $\P$ and $\P^*$.

Recall that the celebrated {\it Logarithmic Barrier Function}
(LBF in short) associated with a
convex optimization problem $\P$ as in (\ref{pbprimal}), or with its dual $\P^*$ in (\ref{dual1}) when 
$g^*$ is explicit, is an important tool
in convex optimization because of its remarkable mathematical properties.
For instance, when the LBF has the {\it self-concordance\footnote{A function $\varphi:\mathcal{D}\to\R$ is called $\kappa$-self-concordant on $\mathcal{D}\subset\R^n$, $\kappa\geq0$, if $\varphi$ is three times continuously differentiable in $\mathcal{D}$, and for all $\x\in\mathcal{D}$ and $h\in\R^n$, one has
\[\vert\nabla^3\varphi(\x)[h,h,h]\vert\leq2\kappa\,\left(h'\nabla^2\varphi(\x)h\right)^{3/2},\]
where $\nabla^3\varphi(\x)[h,h,h]$ is the third differential of $\varphi$ at $\x$
and $h$; see e.g. \cite[p. 52]{hertog}.}} property 
then the associated Logarithmic Barrier algorithm
to solve $\P$ or its dual $\P^*$, runs in time polynomial in the input size of the problem;
see e.g. \cite{guler} and \cite[p. 13, 51 and 60]{hertog}. But the LBF is only one particular choice among many other interior penalty functions! Our main contribution is to 
provide a {\it rationale} behind the LBF as we show
that the LBF can be obtained by approximating
the "$\max$" of a function over a domain by standard $L^p$-norms on the same domain. 
The scalar $1/p$ becomes the parameter of the LBF and nice convergence properties hold when $p\to\infty$. More precisely:

$\bullet$ We first show that the LBF (with parameter $p$) associated with the primal problem $\P$ appears naturally by using the simple and well-known Laplace approximation of a supremum via $L^p$-norms, applied to
the inner infimum in (\ref{min-max1}). It is a bit suprising to obtain an efficient method
in this way. Indeed, the inner infimum in (\ref{min-max1}) (which is exactly equal to $f(\x)$
when $\x$ is an admissible solution of $\P$) is replaced with its "naive" Laplace approximation by $L^p$-norms,
and to the best of our knowledge, the efficiency of this approximation
has not been proved or even tested numerically! 

$\bullet$ Similarly, 
when using the  same Laplace $L^p$-norm approximation technique for the infimum
in the definition (\ref{lfg}) of the conjugate function $g^*$, we obtain a function $\phi_p:\rr^m\to\rr$ which: (a) depends on an 
integer parameter $p$ and (b), is valid on the relative interior $\id$ 
of some domain $\mathcal{D}\subset\rr^m$. Theorem \ref{thmain1} states that
the minimum of $\phi_p$ converges to the minimum of $\P$ as $p\to\infty$.
In doing so for conic 
optimization problems, the set $\mathcal{D}$ is just the 
feasible set of the (known) explicit dual problem $\P^*$, and $\phi_p$ 
is (up to a constant) the LBF with 
parameter $p$, associated with $\P^*$. So again, for conic programs,
the simple Laplace approximation 
of a supremum by $L^p$-norms permits to retrieve the LBF of the dual problem $\P^*$!
Interestingly, Theorem \ref{th-cramer} states that the function 
$\y\mapsto \min_\lambda \phi_p(\lambda;\y)$ is nothing less 
than the {\it Cramer transform} of the Laplace
approximation $\Vert \e^{f}\Vert^p_{L^p(\Omega(\y))}$, 
where $\Omega(\y)\subset\X$ is the feasible set of problem 
$\P$ and $\Vert\cdot\Vert_{L^p}$ is the usual norm associated 
with the Lebesgue space $L^p$.  To the best of our knowledge,
this interpretation of
the Logarithmic Barrier algorithm 
(with parameter $1/p$) for the dual $\P^*$, is new 
(although in the particular context of Linear Programming, 
this result was already alluded to in \cite{or-loga}).

Analogies between the Laplace and Fenchel transforms via exponentials and 
logarithms in the Cramer transform have been already explored in other 
contexts, in order to establish nice parallels between optimization and 
probability via a change of algebra; see e.g. Bacelli et al. \cite{bacel},
Maslov \cite{maslov}, Lasserre \cite{book}, and the many references therein.
In probability, the Cramer transform of a probability measure has also 
been used to provide exact asymptotics of some integrals as well as 
to derive large deviation principles. For a nice survey on this topic the 
interested reader is referred to Piterbarg and Falatov \cite{piter}.

$\bullet$ Finally,  an interesting feature of this Laplace approximation technique
is to provide us with a systematic way to obtain a dual 
problem (\ref{dual1}) in cases when $g^*$ cannot be obtained explicitly
from its definition (\ref{lfg}). Namely, in a number of cases and in contrast with $g^*$,
the function $\phi_p(\lambda;\y)$ obtained by using
the Laplace approximation of the conjugate function $g^*$ by $L^p$-norms, 
can be computed in closed-form explicitly.
Examples of such situations are briefly discussed.
In the general case, $\phi_p$ is of the form 
$h_1(\lambda;\y)+h_2(\lambda;p)$ where: for every $\lambda\in\id$ 
fixed, $h_2(\lambda;p)\to 0$ as $p\to\infty$, and for each 
fixed $p$, the function $\lambda\mapsto h_2(\lambda;p)$ is a barrier 
for the domain $\mathcal{D}$. This yields to consider the optimization 
problem
\[\P^*:\quad \min_\lambda\:\{h_1(\lambda;\y)
\::\:\lambda\in\mathcal{D}\}\]
as a natural dual of $\P$, and for which $\phi_p$ 
is an associated barrier function with parameter $p$. If $g^*$ is concave
then strong duality holds.

\section{Main result}
We need some intermediary helpful results before stating our main result.

\subsection{Some preliminary results}
Let $L^q(\X)$ be the usual Lebesgue space of integrable functions defined 
on a Borel-measurable set $\X\subseteq\rr^n$, and $\Vert h\Vert_{L^q(\X)}$ 
(or sometimes $\Vert h\Vert_q$) be the  associated norm
$$\Vert h\Vert_{L^q(\X)}\,=\,\Vert h\Vert_q\,:=\,
\left(\int_\X\vert h(\x)\vert^q d\x\right)^{1/q}.$$ To make the paper self-contained we
prove the following known result.

\begin{lemma}\label{lemma1}
Let $\X\subseteq\rr^n$ be any Borel-measurable set, and $h\in{}L^q(\X)$ 
for some given $q\geq1$, so that $\|h\|_{L^q(\X)}<\infty$. Then: 
$$\lim_{p\to\infty}\|h\|_{L^p(\X)}\,=\,\Vert{h}
\Vert_\infty\,:=\,\esssup_{\x\in{\X}}|h(\x)|.$$
\end{lemma}

\begin{proof} Notice that $\X$ may be an unbounded set. Suppose that 
$\|h\|_q<\infty$ for some given $q\geq1$, and define $\Lambda$ to be the 
essential suppremum of $|h|$ in $\X$. The result is trivial when 
$\Lambda=0$, so we assume that $\Lambda\in(0,\infty)$. Then
\begin{eqnarray} 
\nonumber
\esssup_{\x\in{\X}}h(\x)&=&\Lambda\;=\;\lim_{p\to\infty}
\big(\|h/\Lambda\|_{L^q(\X)}\big)^{q/p}\Lambda\\
\nonumber&=&\lim_{p\to\infty}\left[\int_\X\Lambda^{p}
\left(\frac{|h(\x)|}{\Lambda}\right)^q\,d\x\right]^{1/p}\\
\label{ineq1}&\geq&\lim_{p\to\infty}\left[\int_\X\vert{h}(\x)
\vert^p\,d\x\right]^{1/p}\,=\,\lim_{p\to\infty}\|h\|_{L^p(\X)}.
\end{eqnarray}
It is also obvious that $\Lambda\geq\lim_p\|h\|_p$ when $\Lambda=\infty$. 
On the other hand, suppose that the essential suppremum $\Lambda$ of 
$|h|$ in $\X$ is finite. Given an arbitrary parameter $\epsilon>0$, there 
exists a bounded subset $B\subset\X$ with positive finite Lebesgue measure 
$\lambda(B)\in(0,\infty)$ such that $|h(\x)|>\Lambda{-}\epsilon$ for every 
$\x\in{B}$. Then
\[\lim_{p\to\infty}\|h\|_{L^p(\X)}\geq\lim_{p\to\infty}\|h\|_{L^p(\B)}\geq
\lim_{p\to\infty}\lambda(B)^{1/p}(\Lambda{-}\epsilon)=\Lambda-\epsilon.\]

Therefore, since $\epsilon$ is arbitrary, combining the 
previous identity with (\ref{ineq1}) yields the desired result 
$\lim_{p\to\infty}\|h\|_{L^p(\X)}=\Lambda$. In the same way, assume that 
the essential supremum of $|h|$ in $\X$ is infinite. Given an arbitrary 
natural number $N\in\N$, there exists a bounded subset $B\subset\X$ with 
positive finite Lebesgue measure $\lambda(B)\in(0,\infty)$ such that 
$|h(\x)|>N$ for every $\x\in{B}$. Then
\[\lim_{p\to\infty}\|h\|_{L^p(\X)}\geq\lim_{p\to\infty}
\|h\|_{L^p(\B)}\geq\lim_{p\to\infty}\lambda(B)^{1/p}N=N.\]

Therefore, since $N$ is arbitrary, combining the 
previous identity with (\ref{ineq1}) yields the desired result 
$\lim_{p\to\infty}\|h\|_{L^p(\X)}=\Lambda=\infty$.
\end{proof}
Next we also need the following intermediate result.
\begin{lemma}
\label{lemma2}
For every $p\in\N$ let $U_p\subset\rr^n$ be some open subset, and let 
$h_p:U_p\to\rr$ be a sequence of functions indexed by the parameter 
$p\in\N$. Suppose that $h_p$ converges pointwise to a function $h$ 
defined on an open subset $U$ of $\rr^n$. Then: 
$$\lim_{p\to\infty}\inf_{\x\in{U_p}}h_p(\x)\,\leq\,\inf_{\x\in{U}}h(\x),$$
provided that the limit in the left side of the equation exists in the 
extended interval $[-\infty,\infty)$. 
\end{lemma}

\begin{proof}
Suppose that the infinimum of $h$ on $U$ is equal to $-\infty$. For every 
$N\in\rr$ there is a point $\x\in{U}$ such that $h(\x)<N$, and so there 
is also an index $p_0$ such that $\x\in{U_p}$ and $h_p(\x)<N$ for every 
$p>p_0$. Hence the infinimum of $h_p$ on $U_p$ is strictly less than $N$, 
and so 
$$\lim_{p\to\infty}\inf_{\x\in{U_p}}h_p(\x)
\,=\,-\infty\,=\,\inf_{\x\in{U}}h(\x),$$
because $N\in\rr$ is arbitrary. On the other hand, assume that the infinimum 
of $h$ on $U$ is equal to $\lambda\in\rr$. For every $\epsilon>0$ there is 
a point $\x\in{U}$ such that $h(\x)<\lambda{+}\epsilon$, and so there is 
also an index $p_0$ such that $\x\in{U_p}$ and $h_p(\x)<\lambda{+}\epsilon$ 
for every $p>p_0$. Since the infinimum of $h_p$ on $U_p$ is strictly less 
than $\lambda{+}\epsilon$ and $\epsilon>0$ is arbitrary, 
$$\lim_{p\to\infty}\inf_{\x\in{U_p}}h_p(\x)
\,\leq\,\lambda\,=\,\inf_{\x\in{U}}h(\x).$$
\end{proof}

\subsection{$L^p$-norm approximations for the primal}

Let us go back to problem $\P$ in (\ref{fenchel1}) where 
$\X\subseteq\R^n$ is a convex cone, and let $\Z:=\R^m_+$.
Let $\X^*\subset\rr^n$ be the dual convex cone associated with $\X$,
and let $\x\mapsto\Delta(\x)$ be the {\it universal logarithmic barrier} function
associated with the convex cone $\X$, that is,
\begin{equation}
\label{univlbf}
\x\mapsto \Delta(\x)\,:=\,\ln\,\left(\int_{\X^*}\e^{-\x'\y} d\y\,\right),\quad\x\in\rix,\end{equation}
where $\rix$ denotes the interior of $\X$. See e.g. G\"uller \cite{guler0} and G\"uler and Tuncel \cite{guler}. Next, let $\mathcal{H}\subset\rr^n$ be the set
\begin{equation}
\label{setH}
\mathcal{H}\,:=\,\{\x\in\rr^n\,:\:\omega(\x)<\y;\: \x\in\rix\}.\end{equation}

Recalling that $\P$ is a maximization problem,
the LBF associated with the (primal) problem $\P$,
and with parameter $p\in\N$, is the function $\psi_p:\mathcal{H}\to\rr$ defined by:
\begin{equation}
\label{lbf-primal}
\x\mapsto\psi_p (\x)\,:=\,f(\x)+\frac{1}{p}\left(-\Delta(\x)+\sum_{j=1}^m\ln (\y-\omega(\x))_j\right).\end{equation}
(In some references like e.g. \cite{hertog}, $p\,\psi_p$ is rather used.)
The LBF in convex programming dates back to 
Frisch \cite{frisch} and became widely known later in Fiacco and McCormick \cite{fiacco}.
For more details and a discussion, see e.g.
den Hertog \cite[Chapter 2]{hertog}.
It is well-known that under some convexity assumptions\footnote{For instance if
the mappings $-f$ and $\omega$ are convex twice continuous differentiable, the interior of
the feasible set is bounded, and
the Hessian $-\nabla^2\psi_p$ is positive definite on its domain; see Den Hertog \cite[page 2]{hertog}. As $p$ varies, the unique maximizer $\x(p)$ of $\psi_p$, called the $p$-{\it center}, lies on the so-called {\it central path}.}, and if $g(\y)<\infty$,
\begin{equation}
\label{lbf-primal1}
g(\y)\,=\,\lim_{p\to\infty}\:\displaystyle\sup_\x\,\{\psi_p(\x)\::\:\x\in\mathcal{H}\,\},
\end{equation}
and the sequence of maximizers $(\x_p)_{p\in\N}\subset\mathcal{H}$ of $\psi_p$ converges to a maximizer of $\P$.

We next provide a simple rationale that explains how the LBF naturally appears to solve problem $\P$. 

\begin{prop}
\label{simple}
With $\mathcal{H}$ as in (\ref{setH}), let $\y\in\R^m$ and $\x\in\mathcal{H}$. Then:
\begin{eqnarray}
\nonumber
0&=&\sup_{(\lambda,\mu)\in\Z\times\X^*}\:\{\lambda' (\omega(\x)-\y)-\mu'\x\,\}\\
\nonumber
&=&\lim_{p\to\infty}\ln\Vert \e^{\lambda'(\omega(\x)-\y)-\mu'\x}\Vert_{L^p(\Z\times\X^*)}\\
&=&\lim_{p\to\infty}\frac{1}{p}\left(\Delta(\x)-\sum_{j=1}^m\ln (\y-\omega(\x))_j\right)
\end{eqnarray}
\end{prop}
\begin{proof}
The first equation is trivial whereas the second one follows from Lemma \ref{lemma1}.
Next,
\begin{eqnarray*}
\ln\Vert \e^{\lambda'(\omega(\x)-\y)-\mu'\x}\Vert_{L^p(\Z\times\X^*)}&=&
\frac{1}{p}
\ln\int_{\Z} \int_{\X^*}\e^{p\lambda'(\omega(\x)-\y)-p\mu'\x}d\mu\,d\lambda\\
&=&\frac{1}{p}\,[\Delta(p\x)-\sum_{j=1}^m\ln (p(\y-\omega(\x)))_j]\\
&=&\frac{1}{p}\,[\Delta(\x)-\sum_{j=1}^m\ln (\y-\omega(\x))_j-(m+n)\ln p]
\end{eqnarray*}
where we have used that for each $p\in\N$, $\Delta(p\x)=\Delta(\x)-n\ln p$ because $\X^*$ is a cone.
\end{proof}
Observe that 
\[\inf_{(\lambda,\mu)\in\Z\times\X^*}\:\{ f(\x)+\lambda' (\y-\omega(\x))+\mu'\x\,\}=
\left\{\begin{array}{cl}
f(\x)&\mbox{if $\x\in\X$ and $\omega(\x)\leq\y$}\\
-\infty &\mbox{otherwise,}\end{array}\right.\]
and so (\ref{min-max1}) can be rewritten
\begin{equation}
\label{rew}
g(\y)=\sup_{\x\in\mathcal{H}}\:\inf_{(\lambda,\mu)\in\Z\times\X^*}\:\{ f(\x)+\lambda' (\y-\omega(\x))+\mu'\x\,\}.
\end{equation}
Using Proposition \ref{simple} and $\lim_{p\to\infty}(m+n)\ln p/p=0$, yields 
\begin{equation}\label{partial1}
g(\y)\,=\,\sup_{\x\in\mathcal{H}}\:\left\{
f(\x)+\lim_{p\to\infty}\frac{1}{p}\left(-\Delta(\x)+\sum_{j=1}^m\ln (\y-\omega(\x))_j\right)
\,\right\}.
\end{equation}
A direct application of Lemma \ref{lemma2} to (\ref{partial1}) yields
\begin{eqnarray*}
g(\y)&\leq&\lim_{p\to\infty}\sup_{\x\in\mathcal{H}}\:\left\{f(\x)+\frac{1}{p}\left(-\Delta(\x)+\sum_{j=1}^m\ln (\y-\omega(\x))_j\right)\right\}\\
&=&\lim_{p\to\infty}\sup_\x\:\{\psi_p(\x)\::\:\x\in\mathcal{H}\},
\end{eqnarray*}
and in fact, (\ref{lbf-primal1}) states that one also has the reverse inequality,

In other words, the LBF $\psi_p$ appears naturally when one approximates
\[\inf_{(\lambda,\mu)\in\Z\times\X^*}\:\{ \lambda' (\y-\omega(\x))+\mu'\x\,\},\]
(whose value is exactly zero if $\omega(\x)\leq \y$ and $\x\in\X$), by the quantity
\[\frac{1}{p}\left(-\Delta(\x)+\sum_{j=1}^m\ln (\y-\omega(\x))_j\right),\]
which comes from the Laplace approximation of a "sup" by $L^p$-norms.

For instance, in Linear programming where $\X=\R^n_+$, 
$\x\mapsto \c'\x$ and $\x\mapsto \omega(\x)=\A\x$ 
for some vector $\c\in\rr^n$ and some matrix $\A\in\rr^{m\times n}$,
\[g(\y)\,=\,\lim_{p\to\infty}\:\displaystyle\sup_\x\,\{ \c'\x+
\frac{1}{p}\left(\sum_{j=1}^m\ln (\y-\A\x)_j+\sum_{i=1}^n\ln (x_i)\right)\,:\,\x\in\mathcal{H}\,\}.\]

\subsection{$L^p$-norm approximations for the dual}

We now use the same approximation technique via $L^p$-norms to either retrieve 
the known dual $\P^*$ when it is explicit, or to provide an explicit dual 
problem $\P^*$  in cases where $g^*$ cannot be obtained explicitly from 
its definition (\ref{lfg}).
Recall that if $g$ is concave, upper semi-continuous, and bounded from above 
by some linear function, then by Legendre-Fenchel duality,
\begin{eqnarray}\label{fenchel2}
g(\y)&=&\inf_\lambda\,\left\{\lambda'\y-
g^*(\lambda)\right\},\quad\hbox{where}\\
\label{fenchel3}g^*(\lambda)&:=&
\inf_\y\,\left\{\lambda'\y-g(\y)\right\}.
\end{eqnarray}
One can express $g^*$ in terms of the definition (\ref{fenchel1}) of $g$ 
and the involved continuous transformations $f$ and $\omega$. Namely, 
\begin{eqnarray}
\nonumber
-g^*(\lambda)&=&\sup_\y\left\{g(\y){-}\lambda'\y\right\}\\
\label{a1}&=&\sup_\y
\sup_{\x\in\X,\atop\omega(\x)\leq\y}\left\{f(\x){-}\lambda'\y\right\}\\
\label{equiv}&=&\left\{\begin{array}{cl}
\displaystyle\sup_{\x\in\X}\left\{f(\x){-}\lambda'
\omega(\x)\right\}&\mbox{if}\;\lambda\geq0,\\ 
+\infty&\mbox{otherwise.}
\end{array}\right.
\end{eqnarray}
Therefore the domain of definition $\mathcal{D}\subset\rr^m$ 
of $g^*$ is given by:
\begin{equation}\label{defd}
\mathcal{D}:=\left\{\lambda\in\rr^m:\lambda\geq0,\;\displaystyle
\sup_{\x\in\X}\left\{f(\x){-}\lambda'\omega(\x)\right\}<\infty\right\},
\end{equation}
with relative interior denoted by $\id$. Observe that $\mathcal{D}$ is convex
because $-g^*$ is convex and proper on $\mathcal{D}$.

\begin{theorem}\label{thmain1}
Let $g$ and $g^*$ be as in (\ref{fenchel1}) and (\ref{fenchel3}), 
respectively. Assume that $g$ is concave, upper semi-continuous, and 
bounded from above by some linear function. 
Suppose that the relative 
interior $\id$ is not empty and for every $\lambda\in\id$ there exists an 
exponent $q\gg1$ such that 
\begin{equation}\label{domain2}
\left\|\e^{f(\x)-\lambda'\omega(\x)}\right\|_{L^q(\X)}\,<\,\infty.
\end{equation}
Then:
\begin{equation}\label{fenchek5}
g(\y)\,=\,\lim_{p\to\infty}\,\inf_{\lambda\in\id}\bigg\{\lambda'\y
+\ln\left\|\e^{f(\x)-\lambda'\omega(\x)}\right\|_{L^p(\X)}\!
-\sum_{j=1}^m\frac{\ln(p\lambda_j)}p\bigg\}
\end{equation}
\end{theorem}

\begin{proof}
In view of (\ref{a1})
\begin{equation}\label{gstar}
-g^*(\lambda)=\left\{\begin{array}{cl}
\ln\bigg[\displaystyle\sup_\y\sup_{\x\in\X\atop \omega(\x)\leq \y}
\left\{\e^{-\lambda'\y+f(\x)}\right\}\bigg]
&\mbox{if}\;\lambda\in\mathcal{D}\\
+\infty&\mbox{otherwise.}
\end{array}\right.\end{equation}

Hypothesis (\ref{domain2}) and Lemma \ref{lemma1} allow 
us to replace the supremum in (\ref{gstar}) by the limit of the 
$L^p$-norms as $p\to\infty$. Namely,
\begin{eqnarray}
\nonumber
-g^*(\lambda)&=&\lim_{p\to\infty}
\ln\left(\int_{\x\in\X}\int_{\omega(\x)\leq \y}
\e^{-p\lambda'\y+pf(\x)}\,d\y\,d\x\right)^{1/p}\\
\nonumber
&=&\lim_{p\to\infty}\left\{
\ln\left(\int_{\x\in\X}
\e^{-p\lambda'\omega(\x)+pf(\x)}\,d\x\right)^{1/p}-\sum_{j=1}^m\frac{\ln (p\lambda_j)}{p}\right\}\\
\label{a2}
&=&\lim_{p\to\infty}
\bigg\{\ln\left\|\e^{f(\x)-\lambda'\omega(\x)}
\right\|_{L^p(\X)}\!-\sum_{j=1}^m\frac{\ln(p\lambda_j)}p\bigg\}.
\end{eqnarray}
Hence from (\ref{a2}), equation (\ref{fenchel2}) 
can be rewritten as follows:
\begin{eqnarray}\nonumber
&&g(\y)\;=\;\inf_{\lambda\in\id}\,\left\{
\lambda'\y-g^*(\lambda)\right\}\;=\\
\nonumber&&=\;\inf_{\lambda\in\id}\,
\lim_{p\to\infty}
\bigg\{\lambda'\y+\ln\left\|\e^{f(\x)-\lambda'\omega(\x)}
\right\|_{L^p(\X)}\!-\sum_{j=1}^m\frac{\ln(p\lambda_j)}p\bigg\}\\
\label{domain3}&&\geq\;\lim_{p\to\infty}\,\inf_{\lambda\in\id}
\bigg\{\lambda'\y+\ln\left\|\e^{f(\x)-\lambda'\omega(\x)}
\right\|_{L^p(\X)}\!-\sum_{j=1}^m\frac{\ln(p\lambda_j)}p\bigg\},
\end{eqnarray}
where we have applied Lemma \ref{lemma2} in order to interchange the 
"$\inf$" and "$\lim$" operators. Notice that the terms between the brackets 
are the functions $h_p(\lambda)$ of Lemma \ref{lemma2}. On the other 
hand, given $\y\in\rr^m$, let 
\[\Theta(\y)\,:=\,\{(\x,\z)\in\X\times\rr^m_+
\,:\,\omega(\x)+\z\leq\y\}\,\subset\,\rr^{n+m},\]
so that whenever $\lambda\in\id$,
\begin{eqnarray*}
\left\|\e^{f(\x)}\right\|_{L^p(\Theta(\y))}&\leq&\left\|\e^{f(\x)
+\lambda'(\y-\omega(\x)-\z)}\right\|_{L^p(\Theta(\y))}\\
&\leq&\e^{\lambda'\y}\left\|\e^{f(\x)-\lambda'\omega(\x)
-\lambda'\z}\right\|_{L^p(\X\times\rr^m_+)}\\
&=&\e^{\lambda'\y}\left\|\e^{f(\x)-\lambda'\omega(\x)}
\right\|_{L^p(\X)}\prod_{j=1}^m(p\lambda_j)^{-1/p}.
\end{eqnarray*}

By hypothesis (\ref{domain2}), given $\lambda\in\id$ fixed, the $L^p$-norm term
in the last above identity is finite for some $p$ large enough. Therefore, 
by definition (\ref{fenchel1}), Lemma~\ref{lemma1}, and the continuity of the logarithm,  one obtains
\begin{eqnarray*}
&&g(\y)\;=\;\ln\,\{\sup_{(\x,\z)\in\Theta(\y)}\e^{f(\x)}\,\}\;=\;
\lim_{p\to\infty}\,\ln\left\|\e^{f(\x)}\right\|_{L^p(\Theta(\y))}\\
&&\leq\;\lim_{p\to\infty}\inf_{\lambda\in\id}\bigg\{\lambda'\y
+\ln\left\|\e^{f(\x)-\lambda'\omega(\x)}\right\|_{L^p(\X)}\!
-\sum_{j=1}^m\frac{\ln(p\lambda_j)}p\bigg\},
\end{eqnarray*}
which combined with (\ref{domain3}) yields the desired result (\ref{fenchek5}).
\end{proof}

\begin{rem}{\rm 
Condition (\ref{domain2}) can be easily checked in particular cases. 
For instance let $\x\mapsto\omega(\x):=\A\x$ for some matrix 
$\A\in\rr^{m\times n}$.
\begin{itemize}
\item If $\X:=\rr^n_+$ and $\x\mapsto{f}(\x):=\c'\x+\ln\x^{\bd}$ for some 
$\c$ and $\bd$ in $\rr^n$ with $\bd\geq0$. The notation $\x^{\bd}$ stands 
for the monomial $\prod_{k=1}^n\x_k^{\bd_k}$. Then $f$ is concave and 
\begin{eqnarray*}
\int_{\rr^n_+}\e^{pf(\x)-p\lambda'\omega(\x)}d\x&=&
\int_{\rr^n_+}\e^{p(\c-\A'\lambda)'\x}\,\x^{p\,\bd}d\x\\
&=&\prod_{k=1}^n\frac{\Gamma(1+p\,\bd_k)}
{(\A_k'\lambda-\c_k)^{1+p\,\bd_k}}\,<\,\infty,\end{eqnarray*}
whenever $\lambda\in\id:=\{\lambda\in\rr^m:\lambda>0,\,\A'\lambda>\c\}$
and $p\in\N$.

\item If $\X=\rr^m$ and $\x\mapsto f(\x):=-\x'\Q\x+\c'\x$ for some 
$\c\in\rr^n$ and a symmetric (strictly) positive definite matrix 
$\Q\in\rr^{n\times n}$, then 
\[\int_{\rr^n}\e^{pf(\x)-p\lambda'\omega(\x)}d\x=\int_{\rr^n}
\e^{-p\x'\Q\x}\,\e^{p(\c-\A'\lambda)'\x}d\x\,<\,\infty,\]
whenever $\lambda\in\id:=\{\lambda\in\X:\lambda>0\}$ and $p\in\N$.
\end{itemize}
}\end{rem}

Consider next the following functions for every $p\in\N$~:
\begin{equation}\label{claim2}
\lambda\mapsto \phi_p(\lambda;\y):=\lambda'\y+\ln
\left\|\e^{f(\x)-\lambda'\omega(\x)}\right\|_{L^p(\X)}
-\sum_{j=1}^m\frac{\ln(p\lambda_j)}p
\end{equation}
defined on some domain of $\rr^m$, and  
\begin{equation}\label{claim1}
\y\mapsto g_p(\y):=\inf_{\lambda}\,\big\{\,
\phi_p(\lambda;\y)\,:\,\lambda\in\id\,\big\}
\end{equation}
defined on $\rr^m$. 
Recall that the {\it Cramer transform} (denoted $\mathcal{C}$) applied 
to an integrable function $u:\rr^m\to\R$, is the Legendre-Fenchel transform 
(denoted $\mathcal{F}$) of the logarithm of the Laplace transform (denoted $\mathcal{L}$) 
of $u$, i.e.,
\[u\mapsto\mathcal{C}(u)\,=\,\mathcal{F}\circ\ln\circ\,\mathcal{L}\,(u).\]

The Cramer transform is natural in the sense that the logarithm of the 
Laplace transform is always a convex function. For our purpose, we 
will consider the concave version of the Fenchel transform
\begin{equation}\label{F-concave}
\hat{u}\,\mapsto\,[\mathcal{F}(\hat{u})](\lambda)
\,=\,\inf_{\y}\{\lambda'\y+\hat{u}(\y)\},
\end{equation}
for $\hat{u}:\rr^m\to\rr$ convex, so that 
$-\hat{u}$ is concave. We claim that:
\begin{theorem}
\label{th-cramer}
The function $\y\mapsto{p}\,g_p(\y)$ defined in 
(\ref{claim1}) is the Cramer transform of the function
\begin{equation}
\label{cramer}
\y\mapsto\tilde{g}_p(\y):=\int_{\Omega(\y)}\!\e^{pf(\x)}
d\x\,=\,\left\Vert\e^f\right\Vert^p_{L^p(\Omega(\y))},
\end{equation}
where $\Omega(\y):=\{\x\in\X:\omega(\x)\leq\y\}\subset\rr^n$.
\end{theorem}

\begin{proof}
The result follows from the definition of 
the Cramer transform $\mathcal{C}$.
\begin{eqnarray*}
\tilde{g}_p\mapsto\mathcal{C}(\tilde{g}_p)&:=&
\mathcal{F}\circ\ln\circ\,\mathcal{L}\,(\tilde{g}_p)\\
\y\mapsto\mathcal{C}(\tilde{g}_p)(\y)&=&\inf_\lambda\big\{
\lambda'\y+[\ln\circ\,\mathcal{L}(\tilde{g}_p)](\lambda)\big\}.
\end{eqnarray*}
Hence
\begin{eqnarray*}
&&[\mathcal{L}(\tilde{g}_p)](p\lambda)\;=\;
\int_{\y\in\rr^m}\e^{-p\lambda' \y}\,\tilde{g}_p(\y)\,d\y\;=\\
&&=\;\int_{\y\in\rr^m}\e^{-p\lambda' \y}\bigg[
\int_{\x\in\X,\,\omega(\x)\leq\y}\e^{pf(\x)}d\x\bigg]d\y\\
&&=\;\int_{\x\in\X}\e^{pf(\x)}\bigg[
\int_{\y\geq\omega(\x)}\e^{-p\lambda' \y}d\y\bigg]d\x\\
&&=\;\bigg[\int_{\x\in\X}\e^{pf(\x)-p\lambda' \omega(\x)}
d\x\bigg]\prod_{j=1}^m\frac1{p\lambda_j}\\
&&=\;\left\Vert\e^{f(\x)-\lambda' \omega(\x)}\right
\Vert^p_{L^p(\X)}\prod_{j=1}^m\frac1{p\lambda_j}.
\end{eqnarray*}
Therefore,
$$[\ln\circ\mathcal{L}(\tilde{g}_p)](p\lambda)=\ln\left\Vert\e^{f(\x)-
\lambda'\omega(\x)}\right\Vert^p_{L^p(\X)}\!-\sum_{j=1}^m\ln(p\lambda_j).$$
On the other hand, recall the definition of 
$g_p(\y)$ given in (\ref{claim1})-(\ref{claim2}), 
$$g_p(\y)=\inf_{\lambda\in\id}\bigg\{\lambda'\y+\ln\left\|\e^{f(\x)-\lambda'
\omega(\x)}\right\|_{L^p(\X)}\!-\sum_{j=1}^m\frac{\ln(p\lambda_j)}p\bigg\}.$$
Thus, with $\mathcal{F}$ as in (\ref{F-concave}) 
and $\mathcal{D}_p:=\{\z:p\,\z\in\mathcal{D}\}$, 
we obtain the desired result~:
\begin{eqnarray*}
p\,g_p(\y)&=&\inf_{\lambda\in\id}\bigg\{p\lambda'\y+\ln\left\|\e^{f(\x)-
\lambda'\omega(\x)}\right\|^p_{L^p(\X)}\!-\sum_{j=1}^m\ln(p\lambda_j)\bigg\}\\
&=&\inf_{\lambda\in\id}\Big\{p\lambda'\y+
[\ln\circ\,\mathcal{L}(\tilde{g}_p)](p\lambda)\Big\}\\
&=&\inf_{\z\in\id_p}\Big\{\z'\y+[\ln\circ
\,\mathcal{L}(\tilde{g}_p)](\z)\Big\}\\
&=&[\mathcal{F}\circ\ln\circ\,\mathcal{L}(\tilde{g}_p)]
(\y)\;=\;[\mathcal{C}(\tilde{g}_p)](\y),
\end{eqnarray*}
\end{proof}

For linear programming, this result was already obtained in \cite{or-loga,book}. 

\begin{ex}(Linear Programming)
\label{ex1}{\rm
In this case set the cone $\X=\rr^n_+$ and the functions $f(\x):=\c'\x$
and $\omega(\x)=\A\x$ for some vector $\c\in\R^n$ and matrix
$\A\in\R^{m\times n}$. We easily have that
$$\left\|\e^{f(\x)-\lambda'\omega(\x)}\right\|^p_{L^p(\X)}=\int_{\X}\e^{p
(\c-\A'\lambda)'\x}d\x=\prod_{k=1}^n\frac1{p\A_k'\lambda-p\,\c_k}$$
for every $p\in\N$ and each $\lambda$ 
in the relative interior $\id$ of the set
\begin{equation}
\mathcal{D}\,=\,\{\lambda\in\R^m\,:
\,\A'\lambda\geq\c,\,\lambda\geq0\}.
\end{equation}
Hence from (\ref{claim2})
\begin{eqnarray*}
\phi_p(\lambda;\y)\;=\;\lambda'\y+\frac1p\ln\left\|
\e^{f(\x)-\lambda'\omega(\x)}\right\|^p_{L^p(\X)}\!
-\sum_{j=1}^m\frac{\ln(p\lambda_j)}p\;=&&\\
=\;\lambda'\y-\sum_{k=1}^n\frac{\ln(\A_k'\lambda-\c_k)}p-
\sum_{j=1}^m\frac{\ln(\lambda_j)}p-\frac{m{+}n}p\ln{p},&&
\end{eqnarray*}

One easily recognizes (up to the constant $(m{+}n)[\ln{p}]/p$)
the LBF with parameter $p$, of the dual problem:
\[\P^*:\quad\min_\lambda\;\{\lambda'\y\,:
\,\A'\lambda\geq\c,\,\lambda\geq0\}.\]
}\end{ex}

\begin{ex}(The general conic problem)
\label{ex2}{\rm
Consider the conic optimization problem
\[\min_\x\;\{\c'\x\,:\,\omega\,\x\leq\y,\,\x\in\X\},\]
for some convex cone $\X\subset\rr^n$, some vector $\c\in\R^n$,
and some linear mapping $\omega:\rr^n\to\rr^m$ with adjoint
mapping $\omega^*:\rr^m\to\rr^n$. We easily have that 
\begin{eqnarray*}
&&\left\|\e^{\c'x-\lambda'\omega\,\x}\right\|^p_{L^p(\X)}\;=\;
\left\|\e^{(\c-\omega^*\lambda)'\x}\right\|^p_{L^p(\X)}\;=\\
&&=\;\int_\X\e^{p(\c-\omega^*\lambda)'\x}d\x\;=\;
p^{-n}\int_\X\e^{(\c-\omega^*\lambda)'\x}d\x,
\end{eqnarray*}
because $\X$ is a cone. Claim (\ref{claim2}) reads
\begin{eqnarray}\nonumber
\phi_p(\lambda;\y)&=&\lambda'\y+\frac1p\ln\left\|
\e^{\c'x-\lambda'\omega\,\x}\right\|^p_{L^p(\X)}\!
-\sum_{j=1}^m\frac{\ln(p\lambda_j)}p\\
\label{phip-conic}&=&\lambda'\y+\frac{\psi(\omega^*\lambda
{-}\c)}p-\sum_{j=1}^m\frac{\ln\lambda_j}p-\frac{m{+}n}p\ln{p},
\end{eqnarray}
where $\psi:\R^n\to\,\R$ is the universal LBF (\ref{univlbf}) associated with the dual cone $\X^*$, and with domain $\id$, where
\begin{equation}
\mathcal{D}\,=\,\{\lambda\in\R^m\,:\,
\omega^*\lambda{-}c\in\X^*,\,\lambda\geq0\}.
\end{equation}
In $\phi_p$ (and up to a constant), one easily recognizes
the LBF with parameter $p$, of the
dual problem:
\[\P^*:\quad\min_\lambda\;\{\lambda'\y\,:\,
\omega^*\lambda{-}c\in\X^*,\,\lambda\geq0\}.\]
}\end{ex}

\begin{ex}(Quadratic programming: non conic formulation)
\label{ex3}{\rm 
Consider symmetric positive semidefinite matrixes 
$\Q_{j}\in\rr^{n\times{n}}$ and vectors $\c_{j}\in\rr^n$ for $j=0,1,...,m$. 
The notation $\Q\succeq0$ (resp. $\Q\succ0$) stands for
$\Q$ is positive semidefinite (resp. strictly positive definite). Let
$\X:=\rr^n$, $f(\x):=-\x'\Q_{0}\x{-}2\c'_{0}\x$ and let
$\omega:\rr^n\to\rr^m$ have entries $\omega_j(\x):=\x'\Q_{j}\x{+}2\c'_{j}\x$
for every $j=1,\ldots,m$. For $\lambda\in\rr^m$ with $\lambda>0$,
define the real symmetric matrix $\Q_{\lambda}\in\rr^{n\times n}$ and vector $\c_{\lambda}\in\rr^n$:
$$\Q_{\lambda}:=\Q_{0}+\sum_{j=1}^m\lambda_j\Q_{j}\quad\hbox
{and}\quad\c_{\lambda}:=\c_{0}+\sum_{j=1}^m\lambda_j\c_{j},$$
so that 
\begin{eqnarray*}
\left\|\e^{f(\x)-\lambda'\omega(\x)}\right\|^p_{L^p(\X)}&=&\int_{\X}
\exp\big({-p}\,\x'\Q_{\lambda}\x{-2p}\,\c'_{\lambda}\x\big)d\x\\
&=&\pi^{n/2}\,\frac{\exp\big(p\,\c'_{\lambda}\Q_{\lambda}^{-1}
\c_{\lambda}\big)}{\sqrt{\det\big(p\,\Q_{\lambda}\big)}}\,<\,\infty,
\end{eqnarray*}
whenever $p\in\N$ and $\Q_{\lambda}\succ0$. Therefore
\begin{eqnarray}\nonumber
\phi_p(\lambda;\y)&=&\lambda'\y+\frac1p\ln\left\|
\e^{\c'x-\lambda'\omega\,\x}\right\|^p_{L^p(\X)}\!
-\sum_{j=1}^m\frac{\ln(p\lambda_j)}p\\
\label{lbf2}&=&\lambda'\y+\c'_{\lambda}\Q_{\lambda}^{-1}
\c_{\lambda}-\frac{\ln\big(\det\Q_{\lambda}\big)}{2p}-\\
\nonumber&&-\sum_{j=1}^m\frac{\ln\lambda_j}p+
n\frac{\ln\pi{-}\ln{p}}{2p}-m\frac{\ln{p}}p,
\end{eqnarray}
on the domain of definition $\id:=\{\lambda:\lambda>0,\,\Q_\lambda\succ0\}$. 
Again, in equation (\ref{lbf2}) one easily recognizes (up to a constant) 
the LBF with parameter $p$, of the dual problem $\P^*$:
\begin{eqnarray*}
\min_{\lambda\geq0,\,\Q_\lambda\succeq0}
&\displaystyle\max_{\x\in \X}&\bigg\{{-}\x'\Q_{0}\x-2\c'_{0}\x-
\sum_{j=1}^m\lambda_j\big(\x'\Q_{j}\x{+}2\c_{j}\x{-}\y_j\big)\bigg\}\\
&=&\min_{\lambda\geq0,\Q_\lambda\succeq0}\left\{\lambda'\y
+\max_{\x\in \X}\big\{-\x'\Q_{\lambda}\x-2\c'_{\lambda}\x)\big\}\right\}\\
&=&\min_{\lambda}\left\{\lambda'\y+\c'_{\lambda}\Q_{\lambda}^{-1}
\c_{\lambda}\,:\,\lambda\geq0,\,\Q_{\lambda}\succeq0\right\},
\end{eqnarray*}
where we have used the fact that $\x^*=\Q_{\lambda}^{-1}\c_{\lambda}\in\rr^n$ 
is the unique optimal solution to the inner maximization problem in 
the second equation above.

If $-\Q_0\succ0$ and $\Q_j\succeq0$, $j=1,\ldots,m$, then
$\id:=\{\lambda:\lambda>0\}$ because $\Q_\lambda\succ0$ whenever $\lambda>0$; in this case $\P$ is a convex optimization problem and there is no duality gap between $\P$ and $\P^*$.

}\end{ex}

\subsection{An explicit dual}
A  minimization problem $\P^*$ in the variables $\lambda\in\mathcal{D}\subset\R^m$ with cost function $\lambda\mapsto h(\lambda)$ is a 
natural dual of $\P$ in (\ref{pbprimal}) if {\it weak duality} holds, that is, if for every feasible solution 
$\lambda\in\mathcal{D}$ of $\P^*$ and every feasible solution $\x\in\R^n$ of $\P$, one has $f(\x)\leq h(\lambda)$. Of course, a highly desirable feature is
that {\it strong duality} holds, that is, the optimal values of $\P$ and $\P^*$ coincide.

In Examples \ref{ex1}, \ref{ex2}, and \ref{ex3},
the function $\phi_p$ defined in (\ref{claim2}) can be 
decomposed into a sum of the form:
\begin{equation}
\label{decomp}
\lambda\mapsto\phi_p(\lambda;\y)\,=\,
h_1(\lambda;\y)+h_2(\lambda;p)\end{equation}
where $h_1$ is independent of the parameter $p$.
Moreover, 
if $h_2(\lambda;p)<\infty$ for some $\lambda>0$ fixed, the term 
$h_2(\lambda;p)$ converges to zero when $p\to\infty$; in addition, for fixed
$p$, $h_2$ is a barrier as $h_2(\lambda;p)\to\infty$ as $\lambda$ approaches the boundary of $\mathcal{D}$.
One may also verify that
\[h_1(\lambda)=\lambda'\y-g^*(\lambda),\quad\forall\lambda\in\mathcal{D},\]
where $g^*$ is the Legendre-Fenchel conjugate in (\ref{fenchel3}).
In fact, the above previous 
decomposition is more general and can be deduced from some simple facts. Recall 
that problem $\P$ is given in (\ref{pbprimal}), and
$$\y\mapsto g(\y)\,:=\,\sup_\x\:\{f(\x):\omega(\x)\leq\y,\,\x\in\X\},$$
for a convex cone $\X\subset\rr^n$ and continuous 
mappings $f$ and $\omega$. Assume that $g$ is concave, upper 
semi-continuous, and bounded from above by some affine function,
so that Legendre-Fenchel duality yields
\begin{eqnarray*}
g(\y)&=&\inf_\lambda\,\left\{\lambda'\y
-g^*(\lambda)\right\},\quad\hbox{where}\\
g^*(\lambda)&:=&\inf_\y\,\left\{\lambda'\y-g(\y)\right\},
\end{eqnarray*}
and where the domain of $g^*$ is the set $\mathcal{D}\subset\rr^m$ 
in (\ref{defd}).
\begin{lemma}\label{lemma3}
Suppose that $\id$ is not empty and for every $\lambda\in\id$ there exists an 
exponent $q\gg1$ such that (\ref{domain2}) holds.
Then:
\begin{equation}
\label{lemma5-1}
\lambda\mapsto h_1(\lambda;\y):=\lim_{p\to\infty}\phi_p(\lambda;\y)\,=\,\lambda'\y-g^*(\lambda),\quad\forall\lambda\in\id.
\end{equation}
\end{lemma}
\begin{proof}
Let $\lambda\in\id$ be fixed. The given hypothesis and Lemma~\ref{lemma1} imply that
\begin{eqnarray*}
&&\lim_{p\to\infty}\ln\left\|\e^{f(\x)-\lambda'\omega(\x)}\right\|_
{L^p(\X)}=\;\sup_{\x\in\X}\left\{f(\x){-}\lambda'\omega(\x)\right\}\;=\\
&&=\;\sup_{\y}\sup_{\x\in\X,\atop\omega(\x)\leq\y}\left\{f(\x){-}\lambda
'\y\right\}\;=\;\sup_\y\left\{g(\y){-}\lambda'\y\right\}\;=\;-g^*(\lambda),
\end{eqnarray*}
and so (\ref{lemma5-1}) follows from the definition (\ref{claim2}) of $\phi_p$.
\end{proof}
As a consequence we obtain:
\begin{cor}
\label{cor-dualpb}
Let $\mathcal{D}$ be as in (\ref{defd}) with $\id\neq\emptyset$, $\phi_p$ as in (\ref{claim2}) and let
$\lambda\mapsto h_1(\lambda;\y)$ be as in (\ref{lemma5-1}).
Then the optimization problem
\begin{equation}
\label{dualpb}
\P^*:\quad \displaystyle\min_\lambda\,\{h_1(\lambda;\y)\::\:\lambda\in\id\:\}.
\end{equation}
is a dual of $\P$. Moreover, if $g$ is concave, upper-semicontinuous and bounded above by some affine function, then strong duality holds.
\end{cor}
\begin{proof}
By Lemma \ref{lemma3}, $h_1(\lambda;\y)=\lambda'\y-g^*(\lambda)$ for all $\lambda\in\id$.
And so if $\min\P^*$ (resp. $\max\P$) denotes the optimal value of $\P^*$ (resp. $\P$), one has
\[\min\P^*=\min_\lambda\,\{\lambda'\y-g^*(\lambda)\::\:\lambda\in\id\}\,\geq\,g(\y)=\max\P,\]
where we have used that $-g^*(\lambda)=\sup_\z\{g(\z)-\lambda'\z\}\geq g(\y)-\lambda'\y$.
Finally, if $g$ is concave, upper semi-continuous and bounded above by some affine function,
then $-g^*$ is convex.
Therefore, as a convex function is continuous 
on its domain $\mathcal{D}$ (which is convex)
\begin{eqnarray*}
\min\P^*&=&\min_\lambda\,\{h_1(\lambda;\y)\::\:\lambda\in\id\:\}\\
&=&\min_\lambda\,\{\lambda'\y-g^*(\lambda)\::\:\lambda\in\id\}\\
&=&\min_\lambda\,\{\lambda'\y-g^*(\lambda)\::\:\lambda\in\mathcal{D}\}\,=\,g(\y),
\end{eqnarray*}
that is, strong duality holds.
\end{proof}
In a number of cases, the $L^p$-norm approximation of $g^*$ can be obtained {\it explicitly} as a function of $\lambda$, whereas
$g^*$ itself cannot be obtained explicitly from (\ref{lfg}). In this situation
one obtains an explicit LBF $\phi_p$ with parameter $p$, for some dual $\P^*$ of $\P$, 
and sometimes 
an explicit dual problem $\P^*$. Indeed if $\phi_p$ is known explicitly, one may sometimes
get its pointwise limit $h_1(\lambda,\y)$ in (\ref{lemma5-1}), in closed form, and so
$\P^*$ is defined explicitly by (\ref{dualpb}). With $p$ fixed, computing $\phi_p(\lambda;\y)$ 
reduces to compute the integral over a convex cone 
of an exponential of some function parametrized by $\lambda$ and $p$. Sometimes this can be done with the help of some known {\it transforms} like e.g. the Laplace or Weierstrass transforms, as illustrated below.

\subsubsection*{Linear mappings and Laplace transform}
Let $\omega:\R^n\to\R^m$ be a linear mapping, with
$\omega(\x)=\A\x$ for some real matrix $\A\in\rr^{m\times n}$, and let $\X=\rr^n_+$. Then
\begin{eqnarray*}
\ln \Vert \e^{f(\x)-\lambda'\omega(\x)}\Vert_{L^p(\X)}&=&
\frac{1}{p}\ln\left(\int_\X\e^{-(p\A'\lambda)'\x}\,e^{pf(\x)}\,d\x\right)\\
&=&\frac{1}{p}\ln\,\left(\mathcal{L}[\e^{pf}](p\A'\lambda)\,\right).
\end{eqnarray*}
That is, the $L^p$-norm approximation is the logarithm of the Laplace transform of the function $\e^f$,
evaluated at the point $p\A'\lambda\in\rr^n$. So if in problem $\P$, the objective function $f$
is such that $\e^f$ has an explicit Laplace transform, then one obtains an explicit expression
for the LBF $\lambda\mapsto \phi_p(\lambda;\y)$ defined in (\ref{claim2}).

For instance if $f(\x)=\c'\x+\ln q(\x)$ for some vector $\c\in\rr^n$ and some polynomial
$q\in\rr[\x]$, positive on the feasible set of $\P$, write
\[q(\x)^p\,=\,\sum_{\alpha\in\N^n}q_{p\alpha}\x^\alpha,\]
for finitely many non zero coefficients $(q_{p\alpha})$, and where the notation 
$\x^\alpha$ stand for the monomial $x_1^{\alpha_1}\cdots x_n^{\alpha_n}$. 
Then 
$\ln\left(\mathcal{L}[\e^{pf}](p\A'\lambda)\right)$ can be computed in closed-form since we have:
\begin{eqnarray*}
\ln\left(\mathcal{L}[\e^{pf}](p\A'\lambda)\right)&=&\ln\left(
\sum_{\alpha\in\N^n}q_{p\alpha}\int_\X\e^{p(c-\A'\lambda)'\x}\,\x^\alpha\,d\x\,\right),\\
&=&-n\ln p+\ln\left(\sum_{\alpha\in\N^n}q_{p\alpha}\frac{\partial^{\vert\alpha\vert}}{\partial 
\x^\alpha}\,\frac{1}{\prod_{i=1}^n(\A'\lambda-c)_i}\right),
\end{eqnarray*}
where $\frac{\partial^{\vert\alpha\vert}}{\partial\x^\alpha}=\prod_{i=1}^n\frac{\partial^{\alpha_i}}{\partial x_i^{\alpha_i}}$.
Of course the above expression can become quite complicated, especially for large values of $p$.
But it is explicit in the variables $(\lambda_i)$. If the function $\x\mapsto \ln q(\x)$ is concave, then 
Corollary \ref{cor-dualpb} applies. On the other hand, to obtain
$g^*$ explicity would require to solve $\A'\lambda-c-\nabla q(\x)/q(\x)=0$ in closed form, 
which is impossible in general.

Similarly if $f$ is linear, i.e. $\x\mapsto f(\x)=\c'\x$ for some vector $\c\in\rr^n$, then
\[\ln \Vert \e^{f(\x)-\lambda'\omega(\x)}\Vert_{L^p(\X)}\,=\,p^{-1}
\ln\left(\mathcal{L}[\e^{-p\lambda'\omega(\x)}](p\c)\right)\]
and so if the function $\x\mapsto \e^{-p\lambda'\omega(\x)}$ has an explicit Laplace transform then so does the $L^p$-norm approximation, and again, $\phi_p$ is obtained in closed form.

\begin{ex}
{\rm As a simple illustrative example, consider the optimization problem:
\begin{equation}
\label{pbex4}
\P:\quad\sup_\x\left\{\c'\x+\sum_{k=1}^n\b_k
\ln\x_k\,:\,\A\x\leq\y,\,\x>0,\,\x\in\rr^n\right\},\end{equation}
for some given matrix $\A\in\R^{m\times n}$ and vectors $\b,\c\in\R^n$ 
and $\y\in\R^m$. We suppose that $\b\geq0$, so that $\c'\x+\ln(\x^\b)$ 
is concave\footnote{The notation $\x^\b$ stands for the monomial $x_1^{b_1}\cdots x_n^{b_n}$.}, and in which case $\P$ is a convex program. 
Notice that with $\X=\rr^n_{++}:=\{\x\in\R^n:\,\x>0\}$,
$$\sup_{\x}\left\{\c'\x+\ln(\x^\b)-\lambda'
\A\x\,:\,\x\in\rr^n,\,\x\in\X\right\}\,<\,\infty$$
whenever $\lambda$ lies in  
$\mathcal{D}=\{\lambda\in\rr^m:\A'\lambda>\c,\,\lambda\geq0\}$. 
\begin{prop}
The function $\phi_p$ in (\ref{claim2}) associated with the optimization problem
(\ref{pbex4}) is given by:
\begin{eqnarray}
\label{propex4-1}
\phi_p(\lambda;\y)&=&\lambda'\y+\sum_{k=1}^n\left[\frac{\ln\Gamma
(1{+}p\,\b_k)}p-\b_k\ln(p\A_k'\lambda{-}p\c_k)\right]\\
\nonumber&&-\sum_{k=1}^n\frac{\ln(\A_k'\lambda{-}\c_k)}p
-\sum_{j=1}^m\frac{\ln\lambda_j}p-\frac{m+n}p\ln{p}.
\end{eqnarray}
and is the LBF with parameter $p$, of the dual problem:
\[\P^*:\quad\inf_\lambda\left\{\lambda'\y-\sum_{k=1}^n
\b_k\ln\left[\frac{\A_k'\lambda{-}\c_k}{\e^{-1}\b_k}\right]
\,:\,\A'\lambda\,>\,\c,\,\lambda\geq0\right\}.\]
In particular, strong duality holds,  i.e., the optimal values of $\P$ and $\P^*$
are equal.
\end{prop}
\begin{proof}
We have
\begin{eqnarray*}
\left\|\e^{\c'\x-\lambda'\A\x}\,\x^{\b}\right\|^p_{L^p(\X)}
&=&\int_{\X}\e^{p(\c-\A'\lambda)'\x}\,\x^{p\,\b}d\x\\
&=&\prod_{k=1}^n\frac{\Gamma(1{+}p\,\b_k)}
{(p\A_k'\lambda{-}p\c_k)^{1+p\,\b_k}}\,<\,\infty,
\end{eqnarray*}
whenever $p\in\N$ and $\lambda>0$ in $\rr^m$ satisfies 
$\A'\lambda>\c$ (and where $\Gamma$ is the usual Gamma function).
Next, $\phi_p$ in (\ref{claim2}) reads
\[\phi_p(\lambda;\y)\,=\,\lambda'\y-\frac1p\ln\left\|
\e^{\c'\x-\lambda'\A\x}\,\x^{\b}\right\|^p_{L^p(\X)}\!
-\sum_{k=1}^n\frac{\ln(p\lambda_k)}p,\]
which is (\ref{propex4-1}).
Next, Stirling's approximation 
$\Gamma(1+t)\approx(t/\e)^t\sqrt{2\pi{t}}$ for real numbers 
$t\gg1$, yields
\begin{eqnarray*}
\lim_{p\to\infty}\,\frac{\ln\Gamma(1{+}p\,\b_k)}p-\b_k\ln{p}&=
&\lim_{p\to\infty}\,\b_k\ln\left[\frac{p\b_k}{\e}\right]-\b_k\ln{p}\\
&=&\b_k\ln(\e^{-1}\b_k).
\end{eqnarray*}
By Lemma~\ref{lemma3},
$$\lambda'\y-g^*(\lambda)=\lim_{p\to\infty}
\phi_p(\lambda;\y)=\lambda'\y-\sum_{k=1}^n\b_k\ln
\left[\frac{\A_k'\lambda{-}\c_k}{\e^{-1}\b_k}\right]\!.$$

And so, by Corollary \ref{cor-dualpb},  the function $\phi_p$ is the LBF with parameter $p$,
of the dual problem
\[\P^*:\quad\inf_\lambda\left\{\lambda'\y-\sum_{k=1}^n
\b_k\ln\left[\frac{\A_k'\lambda{-}\c_k}{\e^{-1}\b_k}\right]
\,:\,\A'\lambda\,>\,\c,\,\lambda\geq0\right\}.\]
In particular, strong duality holds.
\end{proof}
}\end{ex}

\end{document}